\documentclass[12pt]{article}
\usepackage{amssymb}
\usepackage{amsmath}
\usepackage[mathscr]{eucal}

\renewcommand{\phi}{\varphi}

\newcommand{\down}[1]{\stackrel{\raisebox{-0.08cm}{\tiny$\frown$}}{#1}}
\begin{document}
\newcounter{abschnitt}
\newtheorem{koro}{Corollary}[abschnitt]
\newtheorem{defi}[koro]{Definition}
\newtheorem{satz}[koro]{Theorem}
\newtheorem{lem}[koro]{Lemma}
\newtheorem{bem}[koro]{Bemerkungen}
\newtheorem{conj}[koro]{Conjecture}

\newcounter{theoremp}
\newtheorem{theoremp}{Theorem}
\def\thetheoremp{\arabic{abschnitt}.\arabic{theoremp}$\bf_p$}

\newtheorem{korop}[theoremp]{Corollary}
\def\thetheoremp{\arabic{abschnitt}.\arabic{theoremp}$\bf_p$}

\newtheorem{lemp}[theoremp]{Lemma}
\def\thetheoremp{\arabic{abschnitt}.\arabic{theoremp}$\bf_p$}

\newcounter{theoremd}
\newtheorem{theoremd}{Theorem}
\newtheorem{defid}[theoremd]{Definition}
\def\thetheoremd{\arabic{abschnitt}.\arabic{theoremd}$\bf_d$}

\newcounter{saveeqn}
\newcommand{\alpheqn}{\setcounter{saveeqn}{\value{abschnitt}}%
\renewcommand{\theequation}{%
\mbox{\arabic{saveeqn}.\arabic{equation}}}}%
\newcommand{\reseteqn}{\setcounter{equation}{0}%
\renewcommand{\theequation}{\arabic{equation}}}

\begin{center}
{\Large\bf Volume Inequalities and Additive Maps \\
of Convex Bodies}
\end{center}

\vspace{-0.03cm}

\begin{center}
Franz E. Schuster
\end{center}

\begin{center}
Dedicated to Prof.\ Rolf Schneider on the occasion of his 65th
birthday
\end{center}

\vspace{-1.1cm}

\begin{quote}
\footnotesize{ \vskip 1truecm\noindent {\bf Abstract.} Analogs of
the classical inequalities from the Brunn Minkowski Theory for
rotation intertwining additive maps of convex bodies are
developed. We also prove analogs of inequalities from the dual
Brunn Minkowski Theory for intertwining additive maps of star
bodies. These inequalities provide generalizations of results for
projection and intersection bodies. As a corollary we obtain a
new Brunn Minkowski inequality for the volume of polar projection
bodies.

}
\end{quote}

\vspace{1cm}

\centerline{\large{\bf{ \setcounter{abschnitt}{1}
\arabic{abschnitt}. Introduction and Statement of Main Results}}}

\setcounter{theoremp}{2}

\alpheqn

\vspace{0.7cm}\noindent For $n \geq 3$ let $\mathcal{K}^n$ denote
the space of convex bodies (i.e.\ compact, convex sets with
nonempty interior) in $\mathbb{R}^n$ endowed with the Hausdorff
topology. A compact, convex set $K$ is uniquely determined by its
support function $h(K,\cdot)$ on the unit sphere $S^{n-1}$,
defined by $h(K,u)=\max \{u \cdot x: x \in K\}$. If $K \in
\mathcal{K}^n$ contains the origin in its interior, the convex
body $K^*=\{x \in \mathbb{R}^n: x \cdot y \leq 1 \mbox{ for all }
y \in K\}$ is called the polar body of $K$.

The projection body $\Pi K$ of $K \in \mathcal{K}^n$ is the convex
body whose support function is given for $u \in S^{n-1}$ by
\[h(\Pi K,u)=\mbox{vol}_{n-1}(K |u^{\bot}), \]
where $\mbox{vol}_{n-1}$ denotes $(n-1)$-dimensional volume and
$K|u^{\bot}$ is the image of the orthogonal projection of $K$
onto the subspace orthogonal to $u$. Important volume inequalities
for the polars of projection bodies are the Petty projection
inequality \cite{petty72} and the Zhang projection inequality
\cite{zhang91}: Among bodies of given volume the polar projection
bodies have maximal volume precisely for ellipsoids and minimal
volume precisely for simplices. The corresponding results for the
volume of the projection body itself are major open problems in
convex geometry, see \cite{lutwak93b}.  Projection bodies and
their polars have received considerable attention over the last
decades due to their connection to different areas such as
geometric tomography, stereology, combinatorics, computational
and stochastic geometry, see for example
\cite{bourgainlindenstrauss88}, \cite{bolker69},
\cite{gardner95}, \cite{goodeyweil93}, \cite{grinbergzhang99},
\cite{ludwig02}, \cite{ludwig04}, \cite{lutwak85},
\cite{lutwak93}, \cite{schneiderweil83}.

Mixed projection bodies were introduced in the classic volume of
Bonnesen-Fenchel \cite{bonfen34}. They are related to ordinary
projection bodies in the same way that mixed volumes are related
to ordinary volume. In \cite{lutwak85} and \cite{lutwak93} Lutwak
considered the volume of mixed projection bodies and their polars
and established analogs of the classical mixed volume
inequalities.

We will show that the following well known properties of the
projection body operator $\Pi: \mathcal{K}^n \rightarrow
\mathcal{K}^n$ are responsible not only for its behaviour under
Minkowski linear combinations but also for most of the
inequalities established in \cite{lutwak85} and \cite{lutwak93}:
\begin{enumerate}
\item  $\Pi$ is continuous.
\item $\Pi$ is Blaschke Minkowski additive, i.e.\ $\Pi(K \: \# \: L)=\Pi K + \Pi L$ for all $K, L \in
\mathcal{K}^n$.
\item $\Pi$ intertwines rotations, i.e.\ $\Pi(\vartheta K)=\vartheta \Pi K$ for all $K \in \mathcal{K}^n$ and all $\vartheta \in SO(n)$.
\end{enumerate}
Here $\Pi K + \Pi L$ denotes the Minkowski sum of the projection
bodies $\Pi K$ and $\Pi L$ and $K \: \# \: L$ is the Blaschke sum
of the convex bodies $K$ and $L$ (see Section 2). $SO(n)$ is the
group of rotations in $n$ dimensions.

\begin{defi} \label{defbmhomo} A map $\Phi: \mathcal{K}^n \rightarrow
\mathcal{K}^n$ satisfying (a), (b) and (c) is called a Blaschke
Minkowski homomorphism.
\end{defi}

In Section 4 we will see that there are many examples of Blaschke
Minkowski homomorphisms, see also \cite{goodeyweil92},
\cite{schneiderhug02} and \cite{schuster05}. The main purpose of
this article is to extend Lutwak's Brunn Minkowski Theory for
mixed projection operators to general Blaschke Minkowski
homomorphisms. To this end, let in the following $\Phi:
\mathcal{K}^n \rightarrow \mathcal{K}^n$ denote a Blaschke
Minkowski homomorphism.

\begin{satz} \label{mixedbm} There is a continuous operator
\[\Phi: \underbrace{\mathcal{K}^n \times \cdots \times \mathcal{K}^n}_{n-1} \rightarrow \mathcal{K}^n,\]
symmetric in its arguments such that, for $K_1, \ldots, K_m \in
\mathcal{K}^n$ and $\lambda_1, \ldots, \lambda_m \geq 0$,
\begin{equation} \label{expbm}
\Phi(\lambda_1 K_1 + \ldots + \lambda_m K_m)=\sum
\limits_{i_1,\ldots,i_{n-1}} \lambda_{i_1}\cdots \lambda_{i_{n-1}}
\Phi(K_{i_1},\ldots,K_{i_{n-1}}),
\end{equation}
where the sum is with respect to Minkowski addition.
\end{satz}

\noindent Theorem \ref{mixedbm} generalizes the notion of mixed
projection bodies. We will prove a Minkowski inequality for the
volume of mixed Blaschke Minkowski homomorphisms:

\begin{satz} \label{minkbm} If $K, L \in \mathcal{K}^n$, then
\[V(\Phi(K,\ldots,K,L))^{n-1} \geq V(\Phi K)^{n-2}V(\Phi L),   \]
with equality if and only if $K$ and $L$ are homothetic.
\end{satz}

\noindent An Aleksandrov Fenchel type inequality for the volume of
mixed Blaschke Minkowski homomorphisms is provided by:

\begin{satz} \label{alekfenbm} If $K_1, \ldots, K_{n-1} \in \mathcal{K}^n$, then
\[V(\Phi(K_1,\ldots,K_{n-1}))^2\geq V(\Phi(K_1,K_1,K_3,\ldots,K_{n-1}))V(\Phi(K_2,K_2,K_3,\ldots,K_{n-1})).   \]
\end{satz}

\noindent We also prove that the volume of a Blaschke Minkowski
homomorphism satisfies a Brunn Minkowski inequality:

\begin{satz} \label{bmbm} If $K, L \in \mathcal{K}^n$, then
\[V(\Phi(K + L))^{1/n(n-1)} \geq V(\Phi K)^{1/n(n-1)}+V(\Phi L)^{1/n(n-1)},   \]
with equality if and only if $K$ and $L$ are homothetic.
\end{satz}

If we restrict ourselves to even operators, i.e.\ $\Phi K = \Phi
(- K)$ for all $K \in \mathcal{K}^n$, we can also prove volume
inequalities for the polars of mixed Blaschke Minkowski
homomorphisms. In the following we write $\Phi^* K$ for the polar
of $\Phi K$.

\begin{theoremp} \label{minkpbm} If $\Phi$ is even and $K, L \in \mathcal{K}^n$,
then
\[V(\Phi^*(K,\ldots,K,L))^{n-1} \leq V(\Phi^* K)^{n-2}V(\Phi^* L),   \]
with equality if and only if $K$ and $L$ are homothetic.
\end{theoremp}

\noindent This result is again generalized by an Aleksandrov
Fenchel type inequality:

\begin{theoremp} \label{alekfenpbm} If $\Phi$ is even and $K_1, \ldots, K_{n-1} \in \mathcal{K}^n$, then
\[V(\Phi^*(K_1,\ldots,K_{n-1}))^2\leq V(\Phi^*(K_1,K_1,K_3,\ldots,K_{n-1}))V(\Phi^*(K_2,K_2,K_3,\ldots,K_{n-1})).   \]
\end{theoremp}

\noindent The next theorem shows that polars of even Blaschke
Minkowski homomorphisms also satisfy a Brunn Minkowski inequality:

\begin{theoremp} \label{bmpbm} If $\Phi$ is even and $K, L \in \mathcal{K}^n$, then
\[V(\Phi^*(K + L))^{-1/n(n-1)} \geq V(\Phi^* K)^{-1/n(n-1)}+V(\Phi^* L)^{-1/n(n-1)},   \]
with equality if and only if $K$ and $L$ are homothetic.
\end{theoremp}

\noindent Note that the special case $\Phi = \Pi$ of Theorem
\ref{bmpbm} provides a new Brunn Minkowski inequality for the
volume of polar projection bodies:

\begin{korop} If $K, L \in \mathcal{K}^n$, then
\[V(\Pi^*(K + L))^{-1/n(n-1)} \geq V(\Pi^* K)^{-1/n(n-1)}+V(\Pi^* L)^{-1/n(n-1)},   \]
with equality if and only if $K$ and $L$ are homothetic.
\end{korop}

In recent years a dual theory to the Brunn Minkowski Theory of
convex bodies was established. Mixed volumes are replaced by dual
mixed volumes, which were introduced by Lutwak in \cite{lutwak75}.
The natural domain of dual mixed volumes is the space
$\mathcal{S}^n$ of star bodies (i.e.\ compact sets, starshaped
with respect to the origin with continuous radial functions) in
$\mathbb{R}^n$ endowed with the Hausdorff topology. The radial
function $\rho(L,\cdot)$ of a set $L$ starshaped with respect to
the origin is defined on $S^{n-1}$ by $\rho(L,u)=\max\{\lambda
\geq 0: \lambda u \in L\}$.

The intersection body $IL$ of $L \in \mathcal{S}^n$ is the star
body whose radial function is given for $u \in S^{n-1}$ by
\[\rho(IL,u)=\mbox{vol}_{n-1}(L \cap u^{\bot}). \]

Intersection bodies have attracted increased interest in recent
years. They appear already in a paper by Busemann
\cite{busemann53} but were first explicitly defined and named by
Lutwak \cite{lutwak88}. Intersection bodies turned out to be
critical for the solution of the Busemann-Petty problem, see
\cite{gardner94}, \cite{gardner94b}, \cite{gardkolschlum99},
\cite{kalkol05}, \cite{kol98}, \cite{kol00}, \cite{zhang99}. The
fundamental volume inequality for intersection bodies is the
Busemann intersection inequality \cite{busemann53}: Among bodies
of given volume the intersection bodies have maximal volume
precisely for ellipsoids centered in the origin. A corresponding
result for the minimal volume of intersection bodies of a given
volume is another major open problem in convex geometry.

The operator $I: \mathcal{S}^n \rightarrow \mathcal{S}^n$ has the
following well known properties:
\begin{itemize}
\item[$\mbox{(a)}_{\bf d}$] $I$ is continuous.
\item[$\mbox{(b)}_{\bf d}$] $I(K \: \tilde{\#} \: L)=IK \:\tilde{+}\: IL$ for all $K, L \in
\mathcal{S}^n$.
\item[$\mbox{(c)}_{\bf d}$] $I$ intertwines rotations.
\end{itemize}
Here $IK \: \tilde{+} \: IL$ is the radial Minkowski sum of the
intersection bodies $IK$ and $IL$ and $K \: \tilde{\#} \: L$ is
the radial Blaschke sum of the star bodies $K$ and $L$ (see
Section 3).

\begin{defid} \label{defradbmhomo} A map $\Psi: \mathcal{S}^n \rightarrow
\mathcal{S}^n$ is called radial Blaschke Minkowski homomorphism
if it satisfies $\mbox{(a)}_{d}, \mbox{(b)}_d$ and $\mbox{(c)}_d$.
\end{defid}

As Lutwak shows in \cite{lutwak88}, see also \cite{gardner95},
there is a duality between projection and intersection bodies,
that is at present not yet well understood. We will show that
there is a similar duality for general Blaschke Minkowski and
radial Blaschke Minkowski homomorphisms. In analogy to the
inequalities of Theorems \ref{minkbm}, \ref{alekfenbm} and
\ref{bmbm} we will establish dual inequalities for radial Blaschke
Minkowski homomorphisms. To this end, let $\Psi: \mathcal{S}^n
\rightarrow \mathcal{S}^n$ denote a nontrivial radial Blaschke
Minkowski homomorphism, where the operator that maps every star
body to the origin is called the {\it trivial} radial Blaschke
Minkowski homomorphism.

\begin{theoremd} \label{mixdualrbm} There is a continuous operator
\[\Psi: \underbrace{\mathcal{S}^n \times \cdots \times \mathcal{S}^n}_{n-1} \rightarrow \mathcal{S}^n,\]
symmetric in its arguments such that, for $L_1, \ldots, L_m \in
\mathcal{S}^n$ and $\lambda_1, \ldots, \lambda_m \geq 0$,
\[\Psi(\lambda_1 L_1 \: \tilde{+} \: \ldots \: \tilde{+} \: \lambda_m L_m)=\stackrel{\sim}{\sum \limits_{i_1,\ldots,i_{n-1}}} \lambda_{i_1}\cdots \lambda_{i_{n-1}} \Psi(L_{i_1},\ldots,L_{i_{n-1}}),   \]
where the sum is with respect to radial Minkowski addition.
\end{theoremd}

\noindent Mixed intersection bodies were introduced in
\cite{zhang94}. The dual Minkowski inequality for radial Blaschke
Minkowski homomorphisms is:

\begin{theoremd} \label{dualminkrbm} If $K, L \in \mathcal{S}^n$,
then
\[V(\Psi(K,\ldots,K,L))^{n-1} \leq V(\Psi K)^{n-2}V(\Psi L),   \]
with equality if and only if $K$ and $L$ are dilates.
\end{theoremd}

\noindent Theorem \ref{dualminkrbm} is a special case of the dual
Aleksandrov Fenchel inequality for the volume of radial Blaschke
Minkowski homomorphisms:

\begin{theoremd} \label{dualafrbm} If $L_1, \ldots, L_{n-1} \in \mathcal{S}^n$, then
\[V(\Psi(L_1,\ldots,L_{n-1}))^2\leq V(\Psi(L_1,L_1,L_3,\ldots,L_{n-1}))V(\Psi(L_2,L_2,L_3,\ldots,L_{n-1})),   \]
with equality if and only if $L_1$ and $L_2$ are dilates.
\end{theoremd}

\noindent The volume of a radial Blaschke Minkowski homomorphism
satisfies the following dual Brunn Minkowski inequality:

\begin{theoremd} \label{dualbmrbm} If $K, L \in \mathcal{S}^n$, then
\[V(\Psi(K + L))^{1/n(n-1)} \leq V(\Psi K)^{1/n(n-1)}+V(\Psi L)^{1/n(n-1)},   \]
with equality if and only if $K$ and $L$ are dilates.
\end{theoremd}

\noindent Theorems \ref{dualminkrbm}, \ref{dualafrbm} and
\ref{dualbmrbm} for the intersection body operator were recently
established in \cite{lengzhao03} and \cite{lengzhao05}.

\vspace{1cm}

\setcounter{abschnitt}{2}
\centerline{\large{\bf{\arabic{abschnitt}. Mixed Volumes}}}

\reseteqn \alpheqn

\setcounter{koro}{0}

\vspace{0.7cm} \noindent

\noindent We collect in this section the background material and
notation from the Brunn Minkowski Theory that is needed in the
proofs of the main theorems. As a general reference we recommend
the book by Schneider \cite{schneider93}.

The most important algebraic structure on the space
$\mathcal{K}^n$ is Minkowski addition. For $K_1, K_2 \in
\mathcal{K}^n$ and $\lambda_1, \lambda_2 \geq 0$, the support
function of the Minkowski linear combination
$\lambda_1K_1+\lambda_2K_2$ is
\[h(\lambda_1K_1+\lambda_2K_2,\cdot)=\lambda_1h(K_1,\cdot)+\lambda_2h(K_2,\cdot).\]
The volume of a Minkowski linear combination $\lambda_1K_1 +
\ldots + \lambda_m K_m$ of convex bodies $K_1, \ldots, K_m$ can be
expressed as a homogeneous polynomial of degree $n$:
\[V(\lambda_1K_1 + \ldots +\lambda_m K_m)=\sum \limits_{i_1,\ldots, i_n} V(K_{i_1},\ldots,K_{i_n})\lambda_{i_1}\cdots\lambda_{i_n}.   \]
The coefficients $V(K_{i_1},\ldots,K_{i_n})$ are called mixed
volumes of $K_{i_1}, \ldots, K_{i_n}$. These functionals are
nonnegative, symmetric and translation invariant. Moreover, they
are monotone (with respect to set inclusion), multilinear with
respect to Minkowski addition and their diagonal form is ordinary
volume, i.e.\ $V(K,\ldots,K)=V(K)$.

Denote by $V_i(K,L)$ the mixed volume $V(K,\ldots,K,L,\ldots,L)$,
where $K$ appears $n-i$ times and $L$ appears $i$ times. For $0
\leq i \leq n - 1$, write $W_i(K,L)$ for the mixed volume
$V(K,\ldots,K,B,\ldots,B,L)$, where $K$ appears $n-i-1$ times and
the Euclidean unit ball $B$ appears $i$ times. The mixed volume
$W_i(K,K)$ will be written as $W_i(K)$ and is called the $i$th
quermassintegral of $K$. If $\mathbf{C}=(K_1,\ldots,K_i)$, then
$V_i(K,\mathbf{C})$ denotes the mixed volume
$V(K,\ldots,K,K_1,\ldots,K_i)$ with $n-i$ copies of $K$.

Associated with $K_1, \ldots, K_{n-1} \in \mathcal{K}^n$ is a
Borel measure, $S(K_1,\ldots,K_{n-1},\cdot)$, on $S^{n-1}$,
called the mixed surface area measure of $K_1,\ldots, K_{n-1}$.
It is symmetric and has the property that, for each $K \in
\mathcal{K}^n$,

\begin{equation} \label{mixedvolmixedsurf}
V(K,K_1,\ldots,K_{n-1})=\frac{1}{n} \int_{S^{n-1}}
h(K,u)dS(K_1,\ldots,K_{n-1},u) .
\end{equation}

The measures $S_j(K,\cdot):=S(K,\ldots,K,B,\ldots,B,\cdot)$, where
$K$ appears $j$ times and $B$ appears $n-1-j$ times, are called
the surface area measures of order $j$ of $K$. Of particular
importance for our purposes is the surface area measure (of order
$n - 1$) $S_{n-1}(K,\cdot)$ of $K$. It is not concentrated on any
great sphere and has its center of mass in the origin. Conversely,
by Minkowski's existence theorem, every measure in
$\mathcal{M}_+(S^{n-1})$, the space of nonnegative Borel measures
on the sphere with the $\mbox{weak}^*$ topology, with these
properties is the surface area measure of a convex body, unique
up to translation. Hence, if $K_1, K_2 \in \mathcal{K}^n$ and
$\lambda_1, \lambda_2 \geq 0$ (not both 0), then there exists a
convex body $\lambda_1 \cdot K_1 \: \# \: \lambda_2 \cdot K_2$,
unique up to translation, such that

\begin{equation}
S_{n-1}(\lambda_1\cdot K_1 \: \# \: \lambda_2\cdot
K_2,\cdot)=\lambda_1S_{n-1}(K_1,\cdot)+\lambda_2S_{n-1}(K_2,\cdot).
\end{equation}

\noindent This addition and scalar multiplication are called
Blaschke addition and scalar multiplication. For $K \in
\mathcal{K}^n$ and $\lambda \geq 0$, we have $\lambda \cdot K =
\lambda^{1/(n-1)}K.$

Blaschke addition is an additive structure on the space
$[\mathcal{K}^n]$ of translation classes of convex bodies. Thus,
the natural domain of an operator $\Phi$ with the additivity
property $\Phi(K \: \# \: L) = \Phi K + \Phi L$ is the space
$[\mathcal{K}^n]$. In Definition \ref{defbmhomo} the domain of a
Blaschke Minkowski homomorphism is $\mathcal{K}^n$, because we
identify operators on $[\mathcal{K}^n]$ with translation
invariant operators on $\mathcal{K}^n$.

The surface area measure of a Minkowski linear combination of
convex bodies $K_1, \ldots, K_m$ can be expressed as a polynomial
homogeneous of degree $n-1$:

\begin{equation} \label{mixedsurfareameas}
S_{n-1}(\lambda_1K_1+\ldots +\lambda_mK_m,\cdot)=\sum
\limits_{i_1,\ldots,i_{n-1}} \lambda_{i_1}\cdots
\lambda_{i_{n-1}}S(K_{i_1},\ldots,K_{i_{n-1}},\cdot).
\end{equation}

One of the most general and fundamental inequalities for mixed
volumes is the Aleksandrov Fenchel inequality: If $K_1, \ldots,
K_n \in \mathcal{K}^n$ and $1 \leq m \leq n$, then
\begin{equation} \label{alekfen}
V(K_1,\ldots,K_n)^m \geq \prod \limits_{j=1}^m
V(K_j,\ldots,K_j,K_{m+1},\ldots,K_n).
\end{equation}
Unfortunately, the equality conditions of this inequality are, in
general, unknown.

An important special case of inequality (\ref{alekfen}), where the
equality conditions are known, is the Minkowski inequality: If
$K, L \in \mathcal{K}^n$, then
\begin{equation}
V_1(K,L)^n \geq V(K)^{n-1}V(L),
\end{equation}
with equality if and only if $K$ and $L$ are homothetic. In fact,
a more general version of Minkowski's inequality holds: If $0
\leq i \leq n - 2$, then
\begin{equation} \label{genmink}
W_i(K,L)^{n-i} \geq W_i(K)^{n-i-1}W_i(L),
\end{equation}
with equality if and only if $K$ and $L$ are homothetic.

A consequence of the Minkowski inequality is the Brunn Minkowski
inequality: If $K, L \in \mathcal{K}^n$, then
\begin{equation} \label{bmin}
V(K+L)^{1/n} \geq V(K)^{1/n}+V(L)^{1/n},
\end{equation}
with equality if and only if $K$ and $L$ are homothetic. This is
a special case of the more general inequality: If $0 \leq i \leq
n - 2$, then
\begin{equation} \label{quermassbm}
W_i(K+L)^{1/(n-i)} \geq W_i(K)^{1/(n-i)}+W_i(L)^{1/(n-i)},
\end{equation}
with equality if and only if $K$ and $L$ are homothetic.

A further generalization of inequality (\ref{bmin}) is also known
(but without equality conditions): If $0 \leq i \leq n-2$, $K, L,
K_1, \ldots, K_i \in \mathcal{K}^n$ and
$\mathbf{C}=(K_1,...,K_i)$, then
\begin{equation} \label{mostgenbm}
V_i(K+L,\mathbf{C})^{1/(n-i)} \geq
V_i(K,\mathbf{C})^{1/(n-i)}+V_i(L,\mathbf{C})^{1/(n-i)}.
\end{equation}

\vspace{1cm}

\setcounter{abschnitt}{3}
\centerline{\large{\bf{\arabic{abschnitt}. Dual Mixed Volumes}}}

\reseteqn \alpheqn

\setcounter{koro}{0}

\vspace{0.7cm} \noindent

\noindent In this section we summarize results from the dual Brunn
Minkowski Theory, see \cite{lutwak75}. For $L_1, L_2 \in
\mathcal{S}^n$ and $\lambda_1, \lambda_2 \geq 0$, the radial
Minkowski linear combination $\lambda_1
L_1\:\tilde{+}\:\lambda_2L_2$ is the star body defined by
\begin{equation} \label{radminkadd}
\rho(\lambda_1L_1\:\tilde{+}\:\lambda_2L_2,\cdot)=\lambda_1\rho(L_1,\cdot)+\lambda_2\rho(L_2,\cdot).
\end{equation}
The volume of a radial Minkowski linear combination $\lambda_1L_1
\: \tilde{+} \: \ldots \: \tilde{+} \lambda_m L_m$ of star bodies
$L_1, \ldots, L_m$ can be expressed as a homogeneous polynomial
of degree $n$:
\[V(\lambda_1L_1 \: \tilde{+} \: \ldots \: \tilde{+} \: \lambda_m L_m)=\sum \limits_{i_1,\ldots, i_n} \tilde{V}(L_{i_1},\ldots,L_{i_n})\lambda_{i_1}\cdots\lambda_{i_n}.   \]
The coefficients $\tilde{V}(L_{i_1},\ldots,L_{i_n})$ are called
dual mixed volumes of $L_{i_1}, \ldots, L_{i_n}$. They are
nonnegative, symmetric and monotone (with respect to set
inclusion). They are also multilinear with respect to radial
Minkowski addition and $\tilde{V}(L,\ldots,L)=V(L)$. The
following integral representation of dual mixed volumes holds:
\begin{equation} \label{dualmixedvol}
\tilde{V}(L_1,\ldots,L_n)=\frac{1}{n} \int_{S^{n-1}}
\rho(L_1,u)\cdots \rho(L_n,u)du,
\end{equation}
where $du$ is the spherical Lebesgue measure of $S^{n-1}$. The
definitions of $\tilde{V}_i(K,L)$, $\tilde{W}_i(K,L)$, etc. are
analogous to the ones for mixed volumes in Section 2. A slight
extension of the notation $\tilde{V}_i(K,L)$ is for $r \in
\mathbb{R}$
\begin{equation} \label{vminus1}
\tilde{V}_r(K,L)=\frac{1}{n}\int_{S^{n-1}}
\rho^{n-r}(K,u)\rho^r(L,u)du.
\end{equation}

\noindent Obviously we have $\tilde{V}_r(L,L)=V(L)$ for every $r
\in \mathbb{R}$ and every $L \in \mathcal{S}^n$.

\pagebreak

If $\lambda_1, \lambda_2 \geq 0$, then the radial Blaschke linear
combination $\lambda_1 \cdot L_1 \: \tilde{\#} \: \lambda_2 \cdot
L_2$ of the star bodies $L_1$ and $L_2$ is the star body whose
radial function satisfies

\begin{equation}
\rho^{n-1}(\lambda_1\cdot L_1 \: \tilde{\#} \: \lambda_2\cdot
L_2,\cdot)=\lambda_1\rho^{n-1}(L_1,\cdot)+\lambda_2\rho^{n-1}(L_2,\cdot).
\end{equation}

\noindent For $L_1, \ldots, L_m \in \mathcal{S}^{n}$ and
$\lambda_1, \ldots, \lambda_m \geq 0$, the function
$\rho^{n-1}(\lambda_1 L_1 \: \tilde{+} \: \ldots \:\tilde{+}
\lambda_m L_m,\cdot)$  can be expressed as a polynomial
homogeneous of degree $n-1$:

\begin{equation} \label{dualmixedradial}
\rho^{n-1}(\lambda_1L_1 \:\tilde{+}\:\ldots
\:\tilde{+}\:\lambda_mL_m,\cdot)=\sum \limits_{i_1,\ldots,i_{n-1}}
\lambda_{i_1}\cdots \lambda_{i_{n-1}}
\rho(L_{i_1},\cdot)\cdots\rho(L_{i_{n-1}},\cdot).
\end{equation}

The most general inequality for dual mixed volumes is the dual
Aleksandrov Fenchel inequality: If $L_1, \ldots, L_n \in
\mathcal{S}^n$ and $1 \leq m \leq n$, then
\begin{equation} \label{dualalekfen}
\tilde{V}(L_1,\ldots,L_n)^m \leq \prod \limits_{j=1}^m
\tilde{V}(L_j,\ldots,L_j,L_{m+1},\ldots,L_n),
\end{equation}
with equality if and only if $L_1, \ldots, L_m$ are dilates. A
special case of inequality (\ref{dualalekfen}) is the dual
Minkowski inequality: If $K, L \in \mathcal{S}^n$, then
\begin{equation} \label{dualmink}
\tilde{V}_1(K,L) \leq V(K)^{n-1}V(L),
\end{equation}
with equality if and only if $K$ and $L$ are dilates. A more
general version of the dual Minkowski inequality is: If $0 \leq i
\leq n - 2$, then
\begin{equation} \label{dualgenmink}
\tilde{W}_i(K,L)^{n-i} \leq \tilde{W}_i(K)^{n-i-1}\tilde{W}_i(L),
\end{equation}
with equality if and only if $K$ and $L$ are dilates.

We will also need the following Minkowski type inequality: If $K,
L \in \mathcal{S}^n$, then
\begin{equation} \label{minkvminus}
\tilde{V}_{-1}(K,L)^n \geq V(K)^{n+1}V(L)^{-1},
\end{equation}
with equality if and only if $K$ and $L$ are dilates.

A consequence of the dual Minkowski inequality is the dual Brunn
Minkowski inequality: If $K, L \in \mathcal{S}^n$, then
\begin{equation}
V(K \: \tilde{+} \: L)^{1/n} \leq V(K)^{1/n}+V(L)^{1/n},
\end{equation}
with equality if and only if $K$ and $L$ are dilates. Using
Minkowski's integral inequality, this can be further generalized:
If $0 \leq i \leq n - 2$, then
\begin{equation}
\tilde{W}_i(K \: \tilde{+} \: L)^{1/(n-i)} \leq
\tilde{W}_i(K)^{1/(n-i)}+\tilde{W}_i(L)^{1/(n-i)},
\end{equation}
with equality if and only if $K$ and $L$ are dilates. If $0 \leq i
\leq n-2$, $K, L, L_1, \ldots, L_i \in \mathcal{S}^n$ and
$\mathbf{C}=(L_1,...,L_i)$, then
\begin{equation}
\tilde{V}_i(K \: \tilde{+} \: L,\mathbf{C})^{1/(n-i)} \leq
\tilde{V}_i(K,\mathbf{C})^{1/(n-i)}+\tilde{V}_i(L,\mathbf{C})^{1/(n-i)},
\end{equation}
with equality if and only if $K$ and $L$ are dilates.

\pagebreak

\vspace{0.7cm}

\setcounter{abschnitt}{4}
\centerline{\large{\bf{\arabic{abschnitt}. Blaschke Minkowski
homomorphisms}}}

\reseteqn \alpheqn

\setcounter{koro}{0}

\vspace{0.7cm} \noindent In \cite{schneider74},
\cite{schneider74b} Schneider started a systematic investigation
of Minkowski endomorphisms, i.e.\ continuous, rotation
intertwining and Minkowski additive maps of convex bodies. Among
other results he obtained a complete classification of all such
maps in $\mathbb{R}^2$. Kiderlen \cite{kiderlen05} continued these
investigations and extended Schneider's classification to higher
dimensions under a weak monotonicity assumption. He also
classified all Blaschke endomorphisms, i.e.\ continuous, rotation
intertwining and Blaschke additive maps, in arbitrary dimension.
Following the work of Schneider and Kiderlen the author studied
Blaschke Minkowski homomorphisms in \cite{schuster05}. The main
results established there are a representation theorem for
general and a complete classification of all even Blaschke
Minkowski homomorphisms. These results will form the main
ingredients for the proofs of Theorems \ref{mixedbm} to
\ref{bmbm} and Theorems \ref{minkpbm} to \ref{bmpbm}. In order to
state them, we introduce further notation.

$SO(n)$ will be equipped with the invariant probability measure.
As $SO(n)$ is a compact Lie group, the space $\mathcal{M}(SO(n))$
of finite Borel measures on $SO(n)$ with the $\mbox{weak}^*$
topology carries a natural convolution structure. The convolution
$\mu \ast \nu$ of $\mu, \sigma \in \mathcal{M}(SO(n))$ is defined
by
\[\int_{SO(n)} f(\vartheta) d(\mu \ast \sigma)(\vartheta)=\int_{SO(n)} f(\eta \tau)d\mu(\eta)d\sigma(\tau),   \]
for every $f \in \mathcal{C}(SO(n))$, the space of continuous
functions on $SO(n)$ with the uniform topology. By identifying a
continuous function $f$ with the absolute continuous measure with
density $f$, the space $\mathcal{C}(SO(n))$ can be viewed as
subspace of $\mathcal{M}(SO(n))$. Thus the convolution on
$\mathcal{M}(SO(n))$ induces a convolution on
$\mathcal{C}(SO(n))$. Of particular importance for us is the
following Lemma, see \cite{grinbergzhang99}, p.85.

\begin{lem} \label{approxlem} Let $\mu_m, \mu \in \mathcal{M}(SO(n))$,
$m=1,2,\ldots$ and let $f \in \mathcal{C}(SO(n))$. If $\mu_m
\rightarrow \mu$ weakly, then $f \ast \mu_m \rightarrow f \ast
\mu$ and $\mu_m \ast f \rightarrow \mu \ast f$ uniformly.
\end{lem}

Identifying $S^{n-1}$ with the homogeneous space $SO(n)/SO(n-1)$,
where $SO(n-1)$ denotes the group of rotations leaving the point
$\mbox{\raisebox{-0.01cm}{$\down{e}$}}$ (the pole) of $S^{n-1}$
fixed, leads to a one-to-one correspondence of
$\mathcal{C}(S^{n-1})$ and $\mathcal{M}(S^{n-1})$ with right
$SO(n-1)$-invariant functions and measures on $SO(n)$, see
\cite{grinbergzhang99}, \cite{schuster05}. Using this
correspondence, the convolution structure on $\mathcal{M}(SO(n))$
carries over to $\mathcal{M}(S^{n-1})$.

In particular, the convolution $\mu \ast f \in
\mathcal{C}(S^{n-1})$ of a measure $\mu \in \mathcal{M}(SO(n))$
and a function $f \in \mathcal{C}(S^{n-1})$ is defined by
\begin{equation} \label{convvonlinks}
(\mu \ast f)(u)=\int_{SO(n)}f(\vartheta^{-1}u)d\mu(\vartheta).
\end{equation}
If $f=h(K,\cdot)$ is the support function of a compact, convex set
$K$, we have $f(\vartheta^{-1}u)=h(\vartheta K,u)$ for every $u
\in S^{n-1}$. Thus, if $\mu \in \mathcal{M}(SO(n))$ is a
nonnegative measure, $\mu \ast f$ is again the support function of
a compact, convex set which can be interpreted as a weighted
Minkowski rotation mean of the set $K$.

An essential role play convolution operators on
$\mathcal{C}(S^{n-1})$ and $\mathcal{M}(S^{n-1})$, which are
generated by $SO(n-1)$ invariant functions and measures. A
measure $\mu \in \mathcal{M}(S^{n-1})$ is called zonal, if
$\vartheta \mu = \mu$ for every $\vartheta \in SO(n-1)$, where
$\vartheta \mu$ is the image measure under the rotation
$\vartheta$. The set of continuous zonal functions on $S^{n-1}$
will be denoted by
$\mathcal{C}(S^{n-1},\mbox{\raisebox{-0.01cm}{$\down{e}$}})$, the
definition of
$\mathcal{M}(S^{n-1},\mbox{\raisebox{-0.01cm}{$\down{e}$}})$ is
analogous. If $f \in \mathcal{C}(S^{n-1})$, $\mu \in
\mathcal{M}(S^{n-1},\mbox{\raisebox{-0.01cm}{$\down{e}$}})$ and
$\eta \in SO(n)$, we have
\begin{equation} \label{zonalconv}
(f\ast \mu)(\eta \mbox{\raisebox{-0.01cm}{$\down{e}$}})
=\int_{S^{n-1}}f(\eta u)d\mu(u).
\end{equation}
Note that, if $\mu \in
\mathcal{M}(S^{n-1},\mbox{\raisebox{-0.01cm}{$\down{e}$}})$, then,
by (\ref{zonalconv}), for every $f \in \mathcal{C}(S^{n-1})$,
\begin{equation} \label{roteq}
(\vartheta f) \ast \mu = \vartheta(f \ast \mu)
\end{equation}
for every $\vartheta \in SO(n)$. Thus, the spherical convolution
from the right is a rotation intertwining operator on
$\mathcal{C}(S^{n-1})$ and $\mathcal{M}(S^{n-1})$. It is also not
difficult to check from (\ref{zonalconv}) that the convolution of
zonal functions and measures is abelian.

The representation theorem for Blaschke Minkowski homomorphisms
is, see \cite{schuster05}:

\begin{satz} \label{satzbmhomo} If $\Phi: \mathcal{K}^n \rightarrow \mathcal{K}^n$ is a Blaschke
Minkowski homomorphism, then there is a function $g \in
\mathcal{C}(S^{n-1},\mbox{\raisebox{-0.01cm}{$\down{e}$}})$ such
that
\begin{equation} \label{bmhomorep}
h(\Phi K,\cdot)=S_{n-1}(K,\cdot) \ast g.
\end{equation}
The function $g$ is unique up to addition of a function of the
form $u \mapsto x\cdot u, x \in \mathbb{R}^n$.
\end{satz}

We call a compact, convex set $F \subseteq \mathbb{R}^n$ a figure
of revolution if $F$ is invariant under rotations of $SO(n-1)$. A
further investigation of properties of generating functions of
Blaschke Minkowski homomorphisms in \cite{schuster05} led to the
following classification of even Blaschke Minkowski homomorphisms:

\begin{satz} \label{bmsymm} A map $\Phi: \mathcal{K}^n \rightarrow \mathcal{K}^n$
is an even Blaschke Minkowski homomorphism if and only if there is
a centrally symmetric figure of revolution $F \subseteq
\mathbb{R}^n$, which is not a singleton, such that
\[h(\Phi K,\cdot)=S_{n-1}(K,\cdot) \ast h(F,\cdot).\]
The set $F$ is unique up to translations.
\end{satz}

The projection body operator $\Pi: \mathcal{K}^n \rightarrow
\mathcal{K}^n$ is an even Blaschke Minkowski homomorphism. Its
generating figure of revolution is a dilate of the segment
$[-\mbox{\raisebox{-0.01cm}{$\down{e}$}},\mbox{\raisebox{-0.01cm}{$\down{e}$}}]$:
\begin{equation} \label{projbody}
h(\Pi K,\cdot)=\frac{1}{2}S_{n-1}(K,\cdot)\ast
h([-\mbox{\raisebox{-0.01cm}{$\down{e}$}},\mbox{\raisebox{-0.01cm}{$\down{e}$}}],\cdot).
\end{equation}
The operator $\Pi$ maps polytopes to finite Minkowski linear
combinations of rotated and dilated copies of the line segment
$[-\mbox{\raisebox{-0.01cm}{$\down{e}$}},\mbox{\raisebox{-0.01cm}{$\down{e}$}}]$.
A general convex body $K$ is mapped to a limit of such Minkowski
sums of line segments.

Another well known example of an even Blaschke Minkowski
homomorphism is provided by the sine transform of the surface
area measure of a convex body $K$, see \cite{schneiderhug02},
\cite{schneider70}: Define an operator $\Theta: \mathcal{K}^n
\rightarrow \mathcal{K}^n$ by
\[h(\Theta K,\cdot) = S_{n-1}(K,\cdot) \ast h(B \cap
\mbox{\raisebox{-0.01cm}{$\down{e}$}}^{\bot},\cdot).  \] Then
$\Theta$ is an even Blaschke Minkowski homomorphism whose images
are (limits of) Minkowski sums of rotated and dilated copies of
the disc $B \cap \mbox{\raisebox{-0.01cm}{$\down{e}$}}^{\bot}$.
The value $h(\Theta K,u)$ is up to a factor the integrated
surface area of parallel hyperplane sections of $K$ in the
direction $u$.

Every map $\Phi: \mathcal{K}^n \rightarrow \mathcal{K}^n$ of the
form $h(\Phi K,\cdot)=S_{n-1}(K,\cdot) \ast h(L,\cdot)$, for some
figure of revolution $L$, is by (\ref{convvonlinks}) a Blaschke
Minkowski homomorphism, but in general there are generating
functions $g$ of Blaschke Minkowski homomorphisms that are not
support functions. An example of such a map is the (normalized)
second mean section operator $M_2$ introduced in
\cite{goodeyweil92} and further investigated in
\cite{schneiderhug02}: Let $\mathcal{E}_2^n$ be the affine
Grassmanian of two-dimensional planes in $\mathbb{R}^n$ and
$\mu_2$ its motion invariant measure, normalized such that
$\mu_2(\{E \in \mathcal{E}_2^n: E \cap B^n \neq
\varnothing\})=\kappa_{n-2}$. Then
\begin{equation} \label{M2}
h(M_2K,\cdot)=(n-1)\int \limits_{\mathcal{E}_2^n} h(K \cap
E,\cdot)d\mu_2(E)-h(\{z_{n-1}(K)\},\cdot),
\end{equation}
where $z_{n-1}(K)$ is the $(n-1)$st intrinsic moment vector of
$K$, see \cite{schneider93}, p.304.

An immediate consequence of Theorem \ref{satzbmhomo} is Theorem
\ref{mixedbm}. Let $\Phi: \mathcal{K}^n \rightarrow
\mathcal{K}^n$ be a Blaschke Minkowski homomorphism with
generating function $g \in
\mathcal{C}(S^{n-1},\mbox{\raisebox{-0.01cm}{$\down{e}$}})$. If
we define an operator
\[\Phi: \underbrace{\mathcal{K}^n \times \cdots \times \mathcal{K}^n}_{n-1} \rightarrow \mathcal{K}^n,\]
by
\begin{equation} \label{defmixedbm}
h(\Phi(K_1,\ldots,K_{n-1}),\cdot)=S(K_1,\ldots,K_{n-1},\cdot)
\ast g,
\end{equation}
then (\ref{mixedsurfareameas}) and the linearity of convolution
imply (\ref{expbm}). The mixed operator $\Phi$ is well defined
as, by Minkowski's existence theorem, the mixed surface area
measure $S(K_1,\ldots,K_{n-1},\cdot)$ is the surface area measure
(of order $n-1$) of a convex body $[K_1,\ldots,K_{n-1}]$, see
\cite{lutwak86}, and thus
$\Phi(K_1,\ldots,K_{n-1})=\Phi[K_1,\ldots,K_{n-1}].$

By Lemma \ref{approxlem} and the weak continuity of mixed surface
area measures, see \cite{schneider93}, p.276, the mixed operators
defined by (\ref{defmixedbm}) are continuous and symmetric.
Moreover, they have the following properties which are immediate
consequences of the corresponding properties of mixed surface
area measures and the convolution representation
(\ref{defmixedbm}):
\begin{itemize}
\item[(i)] They are multilinear with
respect to Minkowski linear combinations.
\item[(ii)] Their diagonal form reduces to the
Blaschke Minkowski homomorphism:
\[\Phi(K,\ldots,K)=\Phi K.   \]
\item[(iii)] They intertwine simultaneous
rotations, i.e.\ if $\vartheta \in SO(n)$, then
\[\Phi(\vartheta K_1,\ldots,\vartheta K_{n-1})=\vartheta \Phi(K_1,\ldots,K_{n-1}).   \]
\end{itemize}

\noindent For $K, L \in \mathcal{K}^n$, let $\Phi_i(K,L)$ denote
the mixed operator $\Phi(K,\ldots,K,L,\ldots,L)$, with $i$ copies
of $L$ and $n-i-1$ copies of $K$. For the body $\Phi_i(K,B)$ we
simply write $\Phi_iK$.

\pagebreak

\noindent The Steiner point map $s: \mathcal{K}^n \rightarrow
\mathbb{R}^n$ is defined by
\[h(\{s(K)\},\cdot)=nh(K,\cdot) \ast
(\mbox{\raisebox{-0.01cm}{$\down{e}$}} \cdot .\:).\] The map $s$
is the unique vector valued continuous, rigid motion intertwining
and Minkowski additive map on $\mathcal{K}^n$. From the fact that
$S(K_1,\ldots,K_{n-1},\cdot) \ast
(\mbox{\raisebox{-0.01cm}{$\down{e}$}} \cdot .\:)=0$ for $K_1,
\ldots, K_{n-1} \in \mathcal{K}^n$, see \cite{schneider93}, p.281,
we obtain by (\ref{defmixedbm}) and the commutativity of zonal
convolution
\begin{equation*}
h(\{s(\Phi(K_1,\ldots,K_{n-1}))\},\cdot)=nS(K_1,\ldots,K_{n-1},\cdot)
\ast (\mbox{\raisebox{-0.01cm}{$\down{e}$}} \cdot .\:) \ast g = 0.
\end{equation*}
Hence,
\begin{equation} \label{steinerorigin}
s(\Phi(K_1,\ldots,K_{n-1}))=o.
\end{equation}
Since $s(\Phi(K_1,\ldots,K_{n-1})) \in \mbox{int}\,
\Phi(K_1,\ldots,K_{n-1})$, see \cite{schneider93}, p.43, we see
that the convex body $\Phi(K_1,\ldots,K_{n-1})$ contains the
origin in its interior. Thus, the polar body
$\Phi^*(K_1,\ldots,K_{n-1})$, in particular $\Phi^*K$, is well
defined.

For $K \in \mathcal{K}^n$ containing the origin in its interior,
we have the relation $h(K,\cdot)=\rho^{-1}(K^*,\cdot)$. Thus, by
(\ref{defmixedbm}), we obtain for the polar of a mixed Blaschke
Minkowski homomorphism $\Phi$ with generating function $g \in
\mathcal{C}(S^{n-1},\mbox{\raisebox{-0.01cm}{$\down{e}$}})$ the
representation
\begin{equation} \label{rhom1phip}
\rho^{-1}(\Phi^*(K_1,\ldots,K_{n-1}),\cdot)=S(K_1,\ldots,K_{n-1},\cdot)
\ast g.
\end{equation}

\vspace{0.7cm}

\setcounter{abschnitt}{5}
\centerline{\large{\bf{\arabic{abschnitt}. Radial Blaschke
Minkowski homomorphisms}}}

\reseteqn \alpheqn

\setcounter{koro}{0}

\vspace{0.7cm} \noindent In the last section we collected the
representation theorems on Blaschke Minkowski homomorphisms that
are critical in the proofs of Theorems \ref{mixedbm} to
\ref{bmbm} and Theorems \ref{minkpbm} to \ref{bmpbm}. In the
following we will show that there is a corresponding
characterization of radial Blaschke Minkowski homomorphisms that
will be needed to prove the dual Theorems \ref{dualminkrbm},
\ref{dualafrbm} and \ref{dualbmrbm}.

We call a map $\Psi: \mathcal{C}(S^{n-1}) \rightarrow
\mathcal{C}(S^{n-1})$ monotone, if nonnegative functions are
mapped to nonnegative ones. The following theorem is a slight
variation of a result by Dunkl \cite{dunkl66}:

\begin{satz} \label{satzCendos} A map $\Psi: \mathcal{C}(S^{n-1}) \rightarrow \mathcal{C}(S^{n-1})$
is a monotone, linear map that intertwines rotations if and only
if there is a measure $\mu \in
\mathcal{M}_+(S^{n-1},\mbox{\raisebox{-0.01cm}{$\down{e}$}})$ such
that
\begin{equation} \label{Cendos}
\Psi f = f \ast \mu.
\end{equation}
\end{satz}
{\it Proof:} From the definition of spherical convolution and
(\ref{roteq}), it follows that mappings of the form (\ref{Cendos})
have the desired properties.

Conversely, let $\Psi$ be monotone, linear and rotation
intertwining. Consider the map $\psi: \mathcal{C}(S^{n-1})
\rightarrow \mathbb{R},\: f \mapsto \Psi
f(\mbox{\raisebox{-0.01cm}{$\down{e}$}}).$ By the properties of
$\Psi$, the functional $\psi$ is positive and linear on
$\mathcal{C}(S^{n-1})$, thus, by the Riesz representation theorem,
there is a measure $\mu \in \mathcal{M}_+(S^{n-1})$ such that
\[\psi(f)=\int_{S^{n-1}} f(u)d\mu(u).\]
Since $\psi$ is $SO(n-1)$ invariant, the measure $\mu$ is zonal.
Thus, we have for $\eta \in SO(n)$
\[\Psi f(\eta \mbox{\raisebox{-0.01cm}{$\down{e}$}})=
\Psi (\eta^{-1}f)(\mbox{\raisebox{-0.01cm}{$\down{e}$}})=
\psi(\eta^{-1}f)=\int_{S^{n-1}} f(\eta u)d\mu(u).\]
The theorem follows now from (\ref{zonalconv}). \hfill $\blacksquare$\\

\noindent The following consequence of Theorem \ref{satzCendos} is
a dual version of Theorem \ref{satzbmhomo}:

\begin{satz} \label{satzradbmhomo} A map $\Psi: \mathcal{S}^n \rightarrow \mathcal{S}^n$ is a radial Blaschke
Minkowski homomorphism if and only if there is a nonnegative
measure $\mu \in
\mathcal{M}_+(S^{n-1},\mbox{\raisebox{-0.01cm}{$\down{e}$}})$
such that
\begin{equation} \label{radbmhomorep}
\rho(\Psi L,\cdot) = \rho^{n-1}(L,\cdot) \ast \mu.
\end{equation}
\end{satz}
{\it Proof:} From Lemma \ref{approxlem}, (\ref{roteq}) and the
properties of spherical convolution, it is clear that mappings of
the form of (\ref{radbmhomorep}) are radial Blaschke Minkowski
homomorphisms. Thus, we have to show that for every such operator
$\Psi$, there is a measure $\mu \in
\mathcal{M}^+(S^{n-1},\mbox{\raisebox{-0.01cm}{$\down{e}$}})$
such that (\ref{radbmhomorep}) holds.

Since every positive continuous function on $S^{n-1}$ is a radial
function, the vector space
$\{\rho^{n-1}(K,\cdot)-\rho^{n-1}(L,\cdot): K,L \in
\mathcal{S}^n\}$ coincides with $\mathcal{C}(S^{n-1})$. The
operator $\bar{\Psi}: \mathcal{C}(S^{n-1}) \rightarrow
\mathcal{C}(S^{n-1})$ defined by
\[\bar{\Psi}f=\rho(\Psi L_1,\cdot)-\rho(\Psi L_2,\cdot),   \]
where $f = \rho^{n-1}(L_1,\cdot)-\rho^{n-1}(L_2,\cdot)$, is a
linear extension of $\Psi$ to $\mathcal{C}(S^{n-1})$ that
intertwines rotations. Since the cone of radial functions is
invariant under $\bar{\Psi}$, it is also monotone. Hence, by
Theorem \ref{satzCendos}, there is a nonnegative measure $\mu \in
\mathcal{M}_+(S^{n-1},\mbox{\raisebox{-0.01cm}{$\down{e}$}})$ such
that $\bar{\Psi}f=f \ast \mu.$ The statement now follows from
$\bar{\Psi}\rho^{n-1}(L,\cdot)=\rho(\Psi L,\cdot)$.
\hfill $\blacksquare$\\

The generating measure of the intersection body operator $I:
\mathcal{S}^n \rightarrow \mathcal{S}^n$ is the invariant measure
$\mu_{S^{n-2}_0}$ concentrated on $S^{n-2}_0:=S^{n-1} \cap
\mbox{\raisebox{-0.01cm}{$\down{e}$}}^{\bot}$ with total mass
$\kappa_{n-1}$:
\[\rho(IL,\cdot)= \rho^{n-1}(L,\cdot) \ast \mu_{S^{n-2}_0}. \]

Let $\Psi: \mathcal{S}^n \rightarrow \mathcal{S}^n$ be a radial
Blaschke Minkowski homomorphism with generating measure $\mu \in
\mathcal{M}_+(S^{n-1},\mbox{\raisebox{-0.01cm}{$\down{e}$}})$  and
define a mixed operator $\Psi: \mathcal{S}^n \times \cdots \times
\mathcal{S}^n \rightarrow \mathcal{S}^n$ by
\begin{equation} \label{mixedrbm}
\rho(\Psi(L_1,\ldots,L_{n-1}),\cdot)=\rho(L_1,\cdot)\cdots\rho(L_{n-1},\cdot)
\ast \mu.
\end{equation}
The mixed radial Blaschke Minkowski homomorphisms defined in this
way are symmetric and by Lemma \ref{approxlem} continuous.
Moreover, Theorem \ref{mixdualrbm} is a direct consequence of
Theorem \ref{satzradbmhomo} and (\ref{dualmixedradial}). The
properties (ii) and (iii) of mixed Blaschke Minkowski
homomorphisms also hold for mixed radial Blaschke Minkowski
homomorphisms but property (i) has to be replaced by:
\begin{itemize}
\item[$\mbox{(i)}_{\bf d}$] They are multilinear with respect to
radial Minkowski linear combinations.
\end{itemize}
For $K, L \in \mathcal{S}^n$, the definitions of $\Psi_i(K,L)$ and
$\Psi_iK$ are analogous to the ones for mixed Blaschke Minkowski
homomorphisms.

\pagebreak

\vspace{0.7cm}

\setcounter{abschnitt}{6}
\centerline{\large{\bf{\arabic{abschnitt}. Inequalities for
Blaschke Minkowski homomorphisms }}}

\reseteqn \alpheqn

\setcounter{koro}{0}

\vspace{0.7cm} \noindent

\noindent In this section we will prove Theorems \ref{minkbm},
\ref{alekfenbm} and \ref{bmbm} as well as their polar versions.
To this end, let $\Phi: \mathcal{K}^n \rightarrow \mathcal{K}^n$
always denote a Blaschke Minkowski homomorphism with generating
function $g \in
\mathcal{C}(S^{n-1},\mbox{\raisebox{-0.01cm}{$\down{e}$}})$. The
proofs are based on techniques developed by Lutwak in
\cite{lutwak93}.

It will be convenient to introduce the following notation for the
canonical pairing of $f \in \mathcal{C}(S^{n-1})$ and $\mu \in
\mathcal{M}(S^{n-1})$
\begin{equation*}
\langle \mu,f\rangle=\langle f,\mu
\rangle=\int_{S^{n-1}}f(u)d\mu(u).
\end{equation*}

\noindent One very useful tool is the following easy Lemma, see
\cite{schuster05}, p.7:

\begin{lem} \label{selfadlemma}
Let $\mu, \nu \in \mathcal{M}(S^{n-1})$ and $f\in
\mathcal{C}(S^{n-1})$, then
\[\langle \mu \ast \nu ,f \rangle =\langle \mu,f \ast \nu\rangle.\]
\end{lem}

\noindent We summarize geometric consequences of Lemma
\ref{selfadlemma} in the following two Lemmas:

\begin{lem} \label{lemdurch1} If $K_1, \ldots, K_{n-1}, L_1, \ldots,
L_{n-1} \in \mathcal{K}^n$, then
\begin{equation} \label{durch1a}
V(K_1,\ldots,K_{n-1},\Phi(L_1,\ldots,L_{n-1}))=V(L_1,\ldots,L_{n-1},\Phi(K_1,\ldots,K_{n-1})).
\end{equation}
In particular, for $K, L \in \mathcal{K}^n$ and $0 \leq i, j \leq
n-2$,
\begin{equation} \label{durch1b}
W_i(K,\Phi(L_1,\ldots,L_{n-1}))=V(L_1,\ldots,L_{n-1},\Phi_iK)
\end{equation}
and
\begin{equation} \label{durch1c}
W_i(K,\Phi_jL)=W_j(L,\Phi_iK).
\end{equation}
\end{lem}
{\it Proof:} By (\ref{mixedvolmixedsurf}), we have
\[V(K_1,\ldots,K_{n-1},\Phi(L_1,\ldots,L_{n-1}))=
\langle h(\Phi(L_1,\ldots,L_{n-1}),\cdot),S(K_1,\ldots,K_{n-1})
\rangle. \] Hence, identity (\ref{durch1a}) follows from
(\ref{defmixedbm}) and Lemma \ref{selfadlemma}.

For $K_1=\ldots=K_{n-i-1}=K$ and $K_{n-i}=\ldots=K_{n-1}=B$,
identity (\ref{durch1a}) reduces to (\ref{durch1b}). Finally put
$L_1=\ldots=L_{n-j-1}=L$ and $L_{n-j}=\ldots=L_{n-1}=B$ in
(\ref{durch1b}), to obtain identity (\ref{durch1c}). \hfill $\blacksquare$\\

In the next Lemma we summarize further special cases of identity
(\ref{durch1a}). These make use of the fact that the image of a
ball under a Blaschke Minkowski homomorphism is again a ball. To
see this, note that $dS_{n-1}(B,v)=dv$, where $dv$ is the ordinary
spherical Lebesgue measure. Thus, by Theorem \ref{satzbmhomo},
\[h(\Phi B,u)=(S_{n-1}(B,\cdot) \ast g)(u)=\int_{S^{n-1}} g(v)dv=:r_{\Phi}.   \]
So let $r_{\Phi}$ denote the radius of the ball $\Phi B$.

\begin{lem} \label{lemdurch2} If $K_1, \ldots, K_{n-1} \in \mathcal{K}^n$, then
\begin{equation} \label{durch2a}
W_{n-1}(\Phi(K_1,\ldots,K_{n-1}))=r_{\Phi}V(K_1,\ldots,K_{n-1},B).
\end{equation}
In particular, for $K, L \in \mathcal{K}^n$,
\begin{equation} \label{durch2b}
W_{n-1}(\Phi_1(K,L))=r_{\Phi}W_1(K,L),
\end{equation}
and, for $0 \leq i \leq n-2$,
\begin{equation} \label{durch2c}
W_{n-1}(\Phi_iK)=r_{\Phi}W_{i+1}(K).
\end{equation}
\end{lem}

Lemma \ref{lemdurch1} is the critical tool in the proofs of the
inequalities of Theorems \ref{minkbm}, \ref{alekfenbm} and
\ref{bmbm} without the equality conditions, compare Lutwak
\cite{lutwak86}. Lemma \ref{lemdurch2} will be needed to settle
the cases of equality. In fact more general inequalities can be
proved. The following result is a generalization of Theorem
\ref{minkbm}:

\begin{satz} \label{satzgenminkbm} If $K, L \in \mathcal{K}^n$ and $0 \leq i \leq n -
1$, then
\begin{equation} \label{genminkbm}
W_i(\Phi_1(K,L))^{n-1} \geq W_i(\Phi K)^{n-2}W_i(\Phi L),
\end{equation}
with equality if and only if $K$ and $L$ are homothetic.
\end{satz}
{\it Proof:} By (\ref{durch2b}) and (\ref{durch2c}), the case
$i=n-1$ follows from inequality (\ref{genmink}). Let therefore $0
\leq i \leq n - 2$ and $Q \in \mathcal{K}^n$. By (\ref{durch1b})
and (\ref{alekfen}),
\begin{eqnarray*}
W_i(Q,\Phi_1(K,L))^{n-1} & = & V(K,\ldots,K,L,\Phi_i Q)^{n-1}
\geq V_1(K,\Phi_i Q)^{n-2}V_1(L,\Phi_iQ)\\
& = & W_i(Q,\Phi K)^{n-2}W_i(Q,\Phi L).
\end{eqnarray*}
Inequality (\ref{genmink}) implies
\[W_i(Q,\Phi K)^{(n-2)(n-i)}W_i(Q,\Phi L)^{n-i} \geq W_i(Q)^{(n-1)(n-i-1)}W_i(\Phi K)^{n-2}W_i(\Phi L)   \]
and thus,
\begin{equation} \label{equ1}
W_i(Q,\Phi_1(K,L))^{(n-1)(n-i)} \geq
W_i(Q)^{(n-1)(n-i-1)}W_i(\Phi K)^{n-2}W_i(\Phi L),
\end{equation}
with equality if and only if $Q, \Phi K$ and $\Phi L$ are
homothetic. Setting $Q = \Phi_1(K,L)$, we obtain the desired
inequality. If there is equality in (\ref{genminkbm}), we have
equality in (\ref{equ1}). From the fact that the Steiner point of
mixed Blaschke Minkowski homomorphisms is the origin, compare
(\ref{steinerorigin}), it follows that there exist $\lambda_1,
\lambda_2 > 0$ such that
\begin{equation} \label{hom1}
\Phi_1(K,L)=\lambda_1 \Phi K=\lambda_2 \Phi L.
\end{equation}
From the equality in (\ref{genminkbm}), it follows that
\[\lambda_1^{n-2}\lambda_2 = 1.   \]
Moreover, (\ref{durch2b}), (\ref{durch2c}) and (\ref{hom1}) give
\[W_1(K,L)=\lambda_1 W_1(K)=\lambda_2 W_1(L).   \]
Hence, we have
\[W_1(K,L)^{n-1}=W_1(K)^{n-2}W_1(L),  \]
which implies, by (\ref{genmink}), that $K$ and $L$ are
homothetic.
\hfill $\blacksquare$\\

\noindent Of course, Theorem \ref{minkbm} is the special case
$i=0$ of Theorem \ref{satzgenminkbm}.

Much more general then the Minkowski inequality is the
Aleksandrov Fenchel inequality for mixed operators:

\begin{satz} \label{genalekfenbm} If $K_1, \ldots, K_{n-1} \in \mathcal{K}^n$ and $ 1
\leq m \leq n - 1$, then
\[W_i(\Phi(K_1,\ldots,K_{n-1}))^m \geq \prod \limits_{j=1}^m W_i(\Phi(\underbrace{K_j,\ldots,K_j}_m,K_{m+1},\ldots,K_{n-1})).   \]
\end{satz}
{\it Proof:} The case $i=n-1$ reduces by (\ref{durch2a}) to
inequality (\ref{alekfen}). Hence, we can assume $i \leq n - 2$.
From (\ref{durch1b}) and (\ref{alekfen}), it follows that for $Q
\in \mathcal{K}^n$,
\begin{eqnarray*}
W_i(Q,\Phi(K_1,\ldots,K_{n-1}))^m & = &
V(K_1,\ldots,K_{n-1},\Phi_iQ)^m \\
& \geq & \prod \limits_{j=1}^m
V(K_j,\ldots,K_j,K_{m+1},\ldots,K_{n-1},\Phi_iQ) \\
& = & \prod \limits_{j=1}^m
W_i(Q,\Phi(K_j,\ldots,K_j,K_{m+1},\ldots,K_{n-1})).
\end{eqnarray*}
Write $\Phi_{m'}(K_j,\mathbf{C})$ for the mixed operator
$\Phi(K_j,\ldots,K_j,K_{m+1},\ldots,K_{n-1})$. Then, by inequality
(\ref{genmink}), we have
\[W_i(Q,\Phi_{m'}(K_j,\mathbf{C}))^{n-i} \geq
W_i(Q)^{n-i-1}W_i(\Phi_{m'}(K_j,\mathbf{C})).\] Hence, we obtain
\[W_i(Q,\Phi(K_1,\ldots,K_{n-1}))^{m(n-i)} \geq W_i(Q)^{m(n-i-1)}\prod \limits_{j=1}^m
W_i(\Phi_{m'}(K_j,\mathbf{C})).   \] By setting
$Q=\Phi(K_1,\ldots,K_{n-1})$, this becomes the desired inequality.
\hfill $\blacksquare$\\

Theorem \ref{alekfenbm} is the special case $m=2$ and $i=0$ of
Theorem \ref{genalekfenbm}. If we combine the special case
$m=n-2$ of Theorem \ref{genalekfenbm} and Theorem
\ref{satzgenminkbm} we obtain:

\begin{koro} \label{corgenafbm} If $K_1, \ldots, K_{n-1} \in \mathcal{K}^n$ and $0
\leq i \leq n-1$, then
\[W_i(\Phi(K_1,\ldots,K_{n-1}))^{n-1} \geq W_i(\Phi K_1)\cdots W_i(\Phi K_{n-1}),    \]
with equality if and only if the $K_j$ are homothetic.
\end{koro}

The special case $K_1=\ldots=K_{n-1-j}=K$ and
$K_{n-j}=\ldots=K_{n-1}=L$ of Corollary \ref{corgenafbm} leads to
a further generalization of Theorem \ref{minkbm}:

\begin{koro} \label{corgenminkbm} If $K, L \in \mathcal{K}^n$ and $0 \leq i \leq n-1$,
$1 \leq j \leq n-2$, then
\[W_i(\Phi_j(K,L))^{n-1} \geq W_i(\Phi K)^{n-j-1}W_i(\Phi L)^j,   \]
with equality if and only if $K$ and $L$ are homothetic.
\end{koro}

The following theorem provides a general Brunn Minkowski
inequality for the operators $\Phi_j$.

\begin{satz} \label{satzgenbmbm} If $K, L \in \mathcal{K}^n$ and $0 \leq i \leq n-1$,
$0 \leq j \leq n - 3$, then
\begin{equation} \label{genbmbm}
W_i(\Phi_j(K+L))^{1/(n-i)(n-j-1)} \geq
W_i(\Phi_jK)^{1/(n-i)(n-j-1)}+W_i(\Phi_jL)^{1/(n-i)(n-j-1)},
\end{equation}
with equality if and only if $K$ and $L$ are homothetic.
\end{satz}
{\it Proof:} By (\ref{durch1c}) and (\ref{mostgenbm}), we have for
$Q \in \mathcal{K}^n$,
\begin{eqnarray*}
W_i(Q,\Phi_j(K+L))^{1/(n-j-1)} & = & W_j(K+L,\Phi_iQ)^{1/(n-j-1)}
\\
& \geq & W_j(K,\Phi_i Q)^{1/(n-j-1)}+W_j(L,\Phi_i Q)^{1/(n-j-1)}
\\
&=& W_i(Q,\Phi_jK)^{1/(n-j-1)}+W_i(Q,\Phi_jL)^{1/(n-j-1)}.
\end{eqnarray*}
By inequality (\ref{genmink}),
\[W_i(Q,\Phi_jK)^{n-i}\geq W_i(Q)^{n-i-1}W_i(\Phi_jK),   \]
with equality if and only if $Q$ and $\Phi_jK$ are homothetic, and
\[W_i(Q,\Phi_jL)^{n-i}\geq W_i(Q)^{n-i-1}W_i(\Phi_jL),   \]
with equality if and only if $Q$ and $\Phi_jL$ are homothetic.
Thus, we obtain
\begin{eqnarray*}
W_i(Q,\Phi_j(K+L))^{1/(n-j-1)}W_i(Q)^{-(n-i-1)/(n-i)(n-j-1)}
\qquad \qquad \qquad & &
\\ \qquad \qquad \qquad
\geq W_i(\Phi_jK)^{1/(n-i)(n-j-1)}+W_i(\Phi_jL)^{1/(n-i)(n-j-1)},
\end{eqnarray*}
with equality if and only if $Q, \Phi_jK$ and $\Phi_jL$ are
homothetic. If we set $Q=\Phi_j(K+L)$, we obtain (\ref{genbmbm}).
If there is equality in (\ref{genbmbm}), then, by
(\ref{steinerorigin}), there exist $\lambda_1, \lambda_2 > 0$
such that
\begin{equation} \label{hom2}
\Phi_jK=\lambda_1\Phi_j(K+L) \qquad \mbox{and} \qquad
\Phi_jL=\lambda_2 \Phi_j(K+L).
\end{equation}
From equality in (\ref{genbmbm}), it follows that
\[\lambda_1^{1/(n-j-1)}+\lambda_2^{1/(n-j-1)} = 1.   \]
Moreover, (\ref{durch2c}) and (\ref{hom2}) imply
\[W_{j+1}(K)=\lambda_1W_{j+1}(K+L) \qquad \mbox{and} \qquad  W_{j+1}(L)=\lambda_2W_{j+1}(K+L).  \]
Hence, we have
\[W_{j+1}(K+L)^{1/(n-j-1)}=W_{j+1}(K)^{1/(n-j-1)}+W_{j+1}(L)^{1/(n-j-1)},  \]
which implies, by (\ref{quermassbm}), that $K$ and $L$ are
homothetic.
\hfill $\blacksquare$\\

\setcounter{theoremp}{1}

\pagebreak

We turn now to the proofs of Theorems \ref{minkpbm},
\ref{alekfenpbm} and \ref{bmpbm}. To this end, we will {\it
restrict} ourselves to Blaschke Minkowski homomorphisms $\Phi$
with a generating function of the form $g=h(F,\cdot)$, where
$F\subseteq \mathbb{R}^n$ is a figure of revolution which is not
a singleton. Note that, by (\ref{convvonlinks}), every function of
that form is generating function of a Blaschke Minkowski
homomorphism. In particular, by Theorem \ref{bmsymm}, every even
Blaschke Minkowski homomorphism has a generating function of that
type.

We now associate with each such Blaschke Minkowski homomorphism
$\Phi$ a new operator $M_\Phi: \mathcal{S}^n \rightarrow
\mathcal{K}^n$, defined by
\begin{equation} \label{defmphi}
h(M_{\Phi}L,\cdot)=\rho^{n+1}(L,\cdot) \ast h(F,\cdot).
\end{equation}
By (\ref{convvonlinks}), the operator $M_{\Phi}$ is well defined.
Note that $M_{\Phi}$ depends, in contrast to $\Phi$, on the
position of $F$ but that by Theorem \ref{satzbmhomo}, we may
assume that $s(F)=o$. In this \linebreak way we associate to each
Blaschke Minkowski homomorphism a unique operator $M_{\Phi}$.

The next lemma will play the role of Lemma \ref{lemdurch1}.

\begin{lemp} \label{lemdurch3} If $K_1, \ldots, K_{n-1} \in \mathcal{K}^n$ and $L \in \mathcal{S}^n$, then
\begin{equation} \label{durch3a}
\tilde{V}_{-1}(L,\Phi^*(K_1,\ldots,K_{n-1}))=V(K_1,\ldots,K_{n-1},M_{\Phi}L).
\end{equation}
In particular, for $K \in \mathcal{K}^n$,
\begin{equation} \label{durch3b}
\tilde{V}_{-1}(L,\Phi^*_iK)=W_i(K,M_{\Phi}L).
\end{equation}
\end{lemp}
{\it Proof:} By (\ref{vminus1}), we have
\[\tilde{V}_{-1}(K,\Phi^*(K_1,\ldots,K_{n-1}))= \langle \rho^{n+1}(K,\cdot),
\rho^{-1}(\Phi^*(K_1,\ldots,K_{n-1}),\cdot) \rangle.\] Hence,
identity (\ref{durch3a}) follows from (\ref{rhom1phip}) and Lemma
\ref{selfadlemma}. For $K_1=\ldots=K_{n-i-1}=K$ and
$K_{n-i}=\ldots=K_{n-1}=B$,
identity (\ref{durch3a}) reduces to (\ref{durch3b}). \hfill $\blacksquare$\\

\setcounter{theoremp}{3}

We now immediately get the following Minkowski type inequality
for the volume of polar Blaschke Minkowski homomorphisms $\Phi$
with a generating function of the form $g=h(F,\cdot)$. This, in
particular, proves Theorem \ref{minkpbm}:

\begin{theoremp} \label{genminkpbm} If $K, L \in \mathcal{K}^n$,
then
\begin{equation} \label{minkp1}
V(\Phi_1^*(K,L))^{n-1} \leq V(\Phi^* K)^{n-2}V(\Phi^* L),
\end{equation}
with equality if and only if $K$ and $L$ are homothetic.
\end{theoremp}
{\it Proof:} Let $Q \in \mathcal{S}^n$. Then, by (\ref{durch3a})
and (\ref{alekfen}),
\begin{eqnarray*}
\tilde{V}_{-1}(Q,\Phi_1^*(K,L))^{n-1} & = &
V(K,\ldots,K,L,M_{\Phi}Q)^{n-1}
\geq V_1(K,M_{\Phi} Q)^{n-2}V_1(L,M_{\Phi} Q)\\
& = & \tilde{V}_{-1}(Q,\Phi^* K)^{n-2}\tilde{V}_{-1}(Q,\Phi^* L).
\end{eqnarray*}
By inequality (\ref{minkvminus}), we have
\[\tilde{V}_{-1}(Q,\Phi^* K)^{(n-2)n}\tilde{V}_{-1}(Q,\Phi^* L)^n
\geq V(Q)^{(n+1)(n-1)}V(\Phi^*K)^{-(n-2)}V(\Phi^*L)^{-1},   \] and
thus,
\[\tilde{V}_{-1}(Q,\Phi_1^*(K,L))^{(n-1)n}  \geq  V(Q)^{(n+1)(n-1)}V(\Phi^*K)^{-(n-2)}V(\Phi^*L)^{-1},  \]
with equality if and only if $Q, \Phi^* K$ and $\Phi^* L$ are
dilates. Setting $Q = \Phi_1^*(K,L)$, we obtain the desired
inequality. If there is equality in (\ref{minkp1}), then there
exist $\lambda_1, \lambda_2> 0$ such that
\begin{equation} \label{hom16}
\Phi_1^*(K,L)=\lambda_1 \Phi^* K=\lambda_2 \Phi^* L.
\end{equation}
For every convex body $K \in \mathcal{K}^n$ containing the origin
and for every $\lambda > 0$, we have $(\lambda
K)^*=\lambda^{-1}K^*$, and thus
\[\Phi_1(K,L)=\lambda_1^{-1} \Phi K=\lambda_2^{-1} \Phi L.   \]
From the equality in (\ref{minkp1}), it follows that
\[\lambda_1^{-(n-2)}\lambda_2^{-1} = 1.   \]
By (\ref{durch2b}), (\ref{durch2c}) and (\ref{hom16}), we obtain
\[W_1(K,L)=\lambda_1^{-1} W_1(K)=\lambda_2^{-1} W_1(L).   \]
Hence, we have
\[W_1(K,L)^{n-1}=W_1(K)^{n-2}W_1(L),  \]
which implies, by (\ref{genmink}), that $K$ and $L$ are
homothetic.
\hfill $\blacksquare$\\

\begin{theoremp} \label{genalekfenbmp} If $K_1, \ldots, K_{n-1} \in \mathcal{K}^n$ and $ 1
\leq m \leq n - 1$, then
\[V(\Phi^*(K_1,\ldots,K_{n-1}))^m \leq \prod \limits_{j=1}^m V(\Phi^*(\underbrace{K_j,\ldots,K_j}_m,K_{m+1},\ldots,K_{n-1})).   \]
\end{theoremp}
{\it Proof:} From (\ref{durch3a}), it follows that for $Q \in
\mathcal{S}^n$,
\begin{eqnarray*}
\tilde{V}_{-1}(Q,\Phi^*(K_1,\ldots,K_{n-1}))^m & = &
V(K_1,\ldots,K_{n-1},M_{\Phi}Q)^m \\
& \geq & \prod \limits_{j=1}^m
V(K_j,\ldots,K_j,K_{m+1},\ldots,K_{n-1},M_{\Phi}Q) \\
& = & \prod \limits_{j=1}^m
\tilde{V}_{-1}(Q,\Phi^*(K_j,\ldots,K_j,K_{m+1},\ldots,K_{n-1})).
\end{eqnarray*}
Write $\Phi_{m'}^*(K_j,\mathbf{C})$ for the mixed operator
$\Phi^*(K_j,\ldots,K_j,K_{m+1},\ldots,K_{n-1})$. Then, by
inequality (\ref{minkvminus}), we have
\[\tilde{V}_{-1}(Q,\Phi_{m'}^*(K_j,\mathbf{C}))^n \geq
V(Q)^{n+1}V(\Phi_{m'}^*(K_j,\mathbf{C}))^{-1}.\] Hence, we obtain
\[V(Q,\Phi^*(K_1,\ldots,K_{n-1}))^{mn} \geq V(Q)^{m(n+1)}\prod \limits_{j=1}^m
V(\Phi_{m'}^*(K_j,\mathbf{C}))^{-1}.   \] Setting
$Q=\Phi^*(K_1,\ldots,K_{n-1})$, this becomes the desired
inequality.
\hfill $\blacksquare$\\

\pagebreak

Theorem \ref{minkpbm} is the special case $m=2$ of Theorem
\ref{genalekfenbm} for even Blaschke Minkowski homomorphisms.
Combine the special case $m=n-2$ of Theorem \ref{genalekfenbmp}
\linebreak and Theorem \ref{genminkpbm}, to obtain:

\begin{korop} \label{corgenafpbm} If $K_1, \ldots, K_{n-1} \in \mathcal{K}^n$, then
\[V(\Phi^*(K_1,\ldots,K_{n-1}))^{n-1} \leq V(\Phi^* K_1)\cdots V(\Phi^* K_{n-1}),    \]
with equality if and only if the $K_j$ are homothetic.
\end{korop}

The special case, $K_1=\ldots=K_{n-1-j}=K$ and
$K_{n-j}=\ldots=K_{n-1}=L$, of Corollary \ref{corgenafpbm} leads
to a generalization of Theorem \ref{genminkpbm}:

\begin{korop} If $K, L \in \mathcal{K}^n$ and
$1 \leq j \leq n-2$, then
\[V(\Phi_j^*(K,L))^{n-1} \leq V(\Phi^* K)^{n-j-1}V(\Phi^* L)^j,   \]
with equality if and only if $K$ and $L$ are homothetic.
\end{korop}

The last theorem in this section provides a Brunn Minkowski
inequality for the volume of the polar Blaschke Minkowski
homomorphisms under consideration:

\begin{theoremp} \label{satzgenbmpbm} If $K, L \in \mathcal{K}^n$ and $0 \leq j \leq n - 3$, then
\begin{equation} \label{genbmpbm}
V(\Phi_j^*(K+L))^{-1/n(n-j-1)} \geq
V(\Phi_j^*K)^{-1/n(n-j-1)}+V(\Phi_j^*L)^{-1/n(n-j-1)},
\end{equation}
with equality if and only if $K$ and $L$ are homothetic.
\end{theoremp}
{\it Proof:} By (\ref{durch3b}) and (\ref{mostgenbm}), we have for
$Q \in \mathcal{S}^n$,
\begin{eqnarray*}
\tilde{V}_{-1}(Q,\Phi_j^*(K+L))^{1/(n-j-1)} & = &
W_j(K+L,M_{\Phi}Q)^{1/(n-j-1)}
\\
& \geq & W_j(K,M_{\Phi}Q)^{1/(n-j-1)}+W_j(L,M_{\Phi}Q)^{1/(n-j-1)}
\\
&=&
\tilde{V}_{-1}(Q,\Phi_j^*K)^{1/(n-j-1)}+\tilde{V}_{-1}(Q,\Phi_j^*L)^{1/(n-j-1)}.
\end{eqnarray*}
By inequality (\ref{minkvminus}),
\[\tilde{V}_{-1}(Q,\Phi_j^*K)^n \geq V(Q)^{n+1}V(\Phi_j^*K)^{-1},   \]
with equality if and only if $Q$ and $\Phi_j^*K$ are dilates, and
\[\tilde{V}_{-1}(Q,\Phi_j^*L)^n \geq V(Q)^{n+1}V(\Phi_j^*L)^{-1},   \]
with equality if and only if $Q$ and $\Phi_j^*L$ are dilates.
Thus, we obtain
\begin{eqnarray*}
\tilde{V}_{-1}(Q,\Phi_j^*(K+L))^{1/(n-j-1)}V(Q)^{-(n+1)/n(n-j-1)}
\qquad \qquad \qquad & &
\\ \qquad \qquad \qquad
\geq V(\Phi_j^*K)^{-1/n(n-j-1)}+V(\Phi_j^*L)^{-1/n(n-j-1)},
\end{eqnarray*}
with equality if and only if $Q, \Phi_jK$ and $\Phi_jL$ are
dilates. If we set $Q=\Phi_j^*(K+L)$, we obtain (\ref{genbmpbm}).
Suppose equality holds in (\ref{genbmpbm}), then there exist
$\lambda_1, \lambda_2 > 0$ such that
\begin{equation*}
\Phi_j^*K=\lambda_1\Phi_j^*(K+L) \qquad \mbox{and} \qquad
\Phi_j^*L=\lambda_2 \Phi_j^*(K+L),
\end{equation*}
and thus,
\begin{equation} \label{hom17}
\Phi_jK=\lambda_1^{-1}\Phi_j(K+L) \qquad \mbox{and} \qquad
\Phi_jL=\lambda_2^{-1} \Phi_j(K+L).
\end{equation}
From the equality in (\ref{genbmpbm}), it follows that
\[\lambda_1^{-1/(n-j-1)}+\lambda_2^{-1/(n-j-1)} = 1,   \]
and (\ref{durch2c}) and (\ref{hom17}) imply
\[W_{j+1}(K)=\lambda_1^{-1}W_{j+1}(K+L) \qquad \mbox{and} \qquad  W_{j+1}(L)=\lambda_2^{-1}W_{j+1}(K+L).  \]
Hence, we have
\[W_{j+1}(K+L)^{1/(n-j-1)}=W_{j+1}(K)^{1/(n-j-1)}+W_{j+1}(L)^{1/(n-j-1)},  \]
which implies, by (\ref{quermassbm}), that $K$ and $L$ are
homothetic.
\hfill $\blacksquare$\\

\vspace{0.7cm}

\setcounter{abschnitt}{7}
\centerline{\large{\bf{\arabic{abschnitt}. Inequalities for radial
Blaschke Minkowski homomorphisms }}}

\reseteqn \alpheqn

\setcounter{koro}{0}

\vspace{0.7cm}

\noindent The main tools in the proofs of Theorems \ref{minkbm},
\ref{alekfenbm} and \ref{bmbm} are Lemmas \ref{lemdurch1} and
\ref{lemdurch2}. These were immediate consequences of the
convolution representation of Blaschke Minkowski homomorphisms
provided by Theorem \ref{satzbmhomo}. In Section 5, we have shown
that there is a corresponding representation for radial Blaschke
Minkowski homomorphisms, which will now lead to dual versions of
Lemmas \ref{lemdurch1} and \ref{lemdurch2}. In the following let
$\Psi: \mathcal{S}^n \rightarrow \mathcal{S}^n$ denote a
nontrivial radial Blaschke Minkowski homomorphism. In the same
way as Lemmas \ref{lemdurch1} and \ref{lemdurch2} were
consequences of Theorem \ref{satzbmhomo} and Lemma
\ref{selfadlemma}, we obtain from Theorem \ref{satzradbmhomo}:

\begin{lem} \label{lemdurchd1} If $K_1, \ldots, K_{n-1}, L_1, \ldots,
L_{n-1} \in \mathcal{S}^n$, then
\begin{equation} \label{durchd1a}
\tilde{V}(K_1,\ldots,K_{n-1},\Psi(L_1,\ldots,L_{n-1}))=\tilde{V}(L_1,\ldots,L_{n-1},\Psi(K_1,\ldots,K_{n-1})).
\end{equation}
In particular, for $K, L \in \mathcal{S}^n$ and $0 \leq i, j \leq
n-2$,
\begin{equation} \label{durchd1b}
\tilde{W}_i(K,\Psi(L_1,\ldots,L_{n-1}))=\tilde{V}(L_1,\ldots,L_{n-1},\Psi_iK)
\end{equation}
and
\begin{equation} \label{durchd1c}
\tilde{W}_i(K,\Psi_jL)=\tilde{W}_j(L,\Psi_iK).
\end{equation}
\end{lem}

It follows from Theorem \ref{satzradbmhomo} that the image of the
Euclidean unit ball under a radial Blaschke Minkowski
homomorphism $\Psi$ is again a ball. Let $r_{\Psi}$ denote the
radius of this ball. Then the dual version of Lemma
\ref{lemdurch2} is:

\begin{lem} \label{lemdurchd2} If $L_1, \ldots, L_{n-1} \in \mathcal{S}^n$, then
\begin{equation} \label{durchd2a}
\tilde{W}_{n-1}(\Psi(L_1,\ldots,L_{n-1}))=r_{\Psi}\tilde{V}(L_1,\ldots,L_{n-1},B).
\end{equation}
In particular, for $K, L \in \mathcal{S}^n$,
\begin{equation} \label{durchd2b}
\tilde{W}_{n-1}(\Psi_1(K,L))=r_{\Psi}\tilde{W}_1(K,L)
\end{equation}
and, for $0 \leq i \leq n-2$,
\begin{equation} \label{durchd2c}
\tilde{W}_{n-1}(\Psi_iL)=r_{\Psi}\tilde{W}_{i+1}(L).
\end{equation}
\end{lem}

The proofs of Theorems \ref{dualminkrbm}, \ref{dualafrbm} and
\ref{dualbmrbm} are now analogous to the proofs of Theorems
\ref{minkbm}, \ref{alekfenbm} and \ref{bmbm}. We just have to
replace Lemmas \ref{lemdurch1} and \ref{lemdurch2} by Lemmas
\ref{lemdurchd1} and \ref{lemdurchd2}, and use the inequalities
for dual mixed volumes from Section 3 instead of the inequalities
for mixed volumes from Section 2. For this reason we will omit
all the proofs except one in this section:

\begin{satz} \label{genalekfendbm} If $L_1, \ldots, L_{n-1} \in \mathcal{S}^n$ and
$2 \leq m \leq n - 1$, then
\[\tilde{W}_i(\Psi(L_1,\ldots,L_{n-1}))^m \leq \prod \limits_{j=1}^m
\tilde{W}_i(\Psi(\underbrace{L_j,\ldots,L_j}_m,L_{m+1},\ldots,L_{n-1})),
\] with equality if and only if $L_1, \ldots, L_m$ are dilates.
\end{satz}
{\it Proof:} The case $i=n-1$ reduces by (\ref{durch3a}) to
inequality (\ref{dualalekfen}). Hence, assume $i \leq n - 2$. From
(\ref{durchd1a}), it follows that for $Q \in \mathcal{S}^n$,
\begin{eqnarray*}
\tilde{W}_i(Q,\Psi(L_1,\ldots,L_{n-1}))^m & = &
\tilde{V}(L_1,\ldots,L_{n-1},\Psi_iQ)^m \\
& \leq & \prod \limits_{j=1}^m
\tilde{V}(L_j,\ldots,L_j,L_{m+1},\ldots,L_{n-1},\Psi_iQ) \\
& = & \prod \limits_{j=1}^m
\tilde{W}_i(Q,\Psi(L_j,\ldots,L_j,L_{m+1},\ldots,L_{n-1})),
\end{eqnarray*}
with equality if and only if $L_1, \ldots, L_m$ are dilates. Let
$\Psi_{m'}(L_j,\mathbf{C})$ denote the body
$\Psi(L_j,\ldots,L_j,L_{m+1},\ldots,L_{n-1})$. Then, by inequality
(\ref{dualgenmink}), we have
\[\tilde{W}_i(Q,\Psi_{m'}(L_j,\mathbf{C}))^{n-i} \leq
\tilde{W}_i(Q)^{n-i-1}\tilde{W}_i(\Psi_{m'}(L_j,\mathbf{C})),\]
with equality if and only if $Q$ and $\Psi_{m'}(L_j,\mathbf{C})$
are dilates. Hence,
\[\tilde{W}_i(Q,\Psi^*(L_1,\ldots,L_{n-1}))^{m(n-i)} \leq \tilde{W}_i(Q)^{m(n-i-1)}\prod \limits_{j=1}^m
\tilde{W}_i(\Psi_{m'}(L_j,\mathbf{C})).   \] By setting
$Q=\Psi(L_1,\ldots,L_{n-1})$, the statement follows.
\hfill $\blacksquare$\\

Theorem \ref{dualminkrbm} and \ref{dualafrbm} are now just special
cases of Theorem \ref{genalekfendbm}. Further consequences are the
dual versions of Corollaries \ref{corgenafbm} and
\ref{corgenminkbm}:

\begin{koro} If $L_1, \ldots, L_{n-1} \in \mathcal{S}^n$ and $0
\leq i \leq n-1$, then
\[\tilde{W}_i(\Psi(L_1,\ldots,L_{n-1}))^{n-1} \leq \tilde{W}_i(\Psi L_1)\cdots \tilde{W}_i(\Psi L_{n-1}),    \]
with equality if and only if the $L_j$ are dilates.
\end{koro}

\begin{koro} If $K, L \in \mathcal{S}^n$ and $0 \leq i \leq n-1$,
$1 \leq j \leq n-2$, then
\[\tilde{W}_i(\Psi_j(K,L))^{n-1} \leq \tilde{W}_i(\Psi K)^{n-j-1}\tilde{W}_i(\Psi L)^j,   \]
with equality if and only if $K$ and $L$ are dilates.
\end{koro}

The dual counterpart of Theorem \ref{satzgenbmbm} is:

\begin{satz} If $K, L \in \mathcal{S}^n$ and $0 \leq i \leq n-1$,
$0 \leq j \leq n - 3$, then
\begin{equation}
\tilde{W}_i(\Psi_j(K+L))^{1/(n-i)(n-j-1)} \leq
\tilde{W}_i(\Psi_jK)^{1/(n-i)(n-j-1)}+\tilde{W}_i(\Phi_jL)^{1/(n-i)(n-j-1)},
\end{equation}
with equality if and only if $K$ and $L$ are dilates.
\end{satz}

\pagebreak

\vspace{0.7cm}

\setcounter{abschnitt}{8}
\centerline{\large{\bf{\arabic{abschnitt}. Final remarks }}}

\reseteqn \alpheqn

\setcounter{koro}{0}

\vspace{0.7cm} \noindent

\noindent In Theorems \ref{genminkpbm}, \ref{genalekfenbmp} and
\ref{satzgenbmpbm}, we restrict ourselves to Blaschke Minkowski
homomorphisms $\Phi$ with generating functions $g$ that are
support functions. We did this to ensure that star bodies are
mapped to convex bodies by the operators $M_{\Phi}$ defined in
(\ref{defmphi}). An example of a Blaschke Minkowski homomorphism
whose generating function is not a support function is the second
mean section operator $M_2$, see (\ref{M2}). A natural question
is whether Theorems \ref{genminkpbm}, \ref{genalekfenbmp} and
\ref{satzgenbmpbm} hold for general Blaschke Minkowski
homomorphisms.

If $\Phi$ is the projection body operator, the map $M_{\Phi}$
becomes a multiple of the moment body operator which is (up to
volume normalization) the well known centroid body operator
$\Gamma: \mathcal{S}^n \rightarrow \mathcal{K}^n$. Centroid
bodies were defined and investigated by Petty \cite{petty61}.
They have proven to be an important tool in establishing
fundamental affine isoperimetric inequalities, see
\cite{gardner95}, \cite{lutwak90}, \cite{milmanpajor89},
\cite{petty72}. The Busemann-Petty centroid inequality, for
example, states that
\begin{equation} \label{busepet}
V(\Gamma K) \geq \left( \frac{2\kappa_{n-1}}{(n+1)\kappa_n}
\right)^nV(K),
\end{equation}
where $\kappa_n$ is the volume of the Euclidean unit ball in $n$
dimensions. Inequality (\ref{busepet}) is critical for the proof
of Petty's projection inequality
\begin{equation} \label{petproj}
V(K)^{n-1}V(\Pi^*K) \leq \left( \frac{\kappa_n}{\kappa_{n-1}}
\right)^n.
\end{equation}
It is the author's belief that an inequality corresponding to
(\ref{busepet}) holds for all operators $M_\Phi$. This would
immediately provide a generalization of Petty's inequality to
general Blaschke Minkowski homomorphisms and would show that the
affine invariant inequality (\ref{petproj}) holds in a more
general setting.

\vspace{0.7cm}
\medskip\noindent{\bf Acknowledgements.} The work of
the author was supported by the Austrian Science Fund (FWF),
within the scope of the project "Affinely associated bodies",
Project Number: P16547-N12 and the project "Phenomena in high
dimensions" of the European Community, Contract Number:
MRTN-CT-2004-511953. The author is obliged to Monika Ludwig for
her helpful remarks.

\hfill\parbox[t]{6truecm}{ Forschungsgruppe \hfill\par Konvexe
und Diskrete Geometrie \hfill\par Technische Universit\"at
Wien\hfill\par Wiedner Hauptstra\ss e 8--10/1046\hfill\par A--1040
Vienna, Austria\hfill\par franz.schuster@tuwien.ac.at\hfill}

\end{document}